\documentclass[letterpaper,11pt]{article}
\usepackage{amssymb}
\usepackage{amsmath}
\usepackage{amsthm}
\usepackage[import]{xypic}
\usepackage{graphicx}

\def\g{\mathfrak g}      
\def\h{\mathfrak h}      
\def\u{\mathfrak u}      
\def\t{\mathfrak t}      
\def\k{\mathfrak k}      
\def\n{\mathfrak n}      
\def\p{\mathfrak p}
\def\d{\mathfrak d}      
\def\a{\mathfrak a}
\def\l{\mathfrak l}

\def\H{\mathcal H}

\def\w{\mathbf w}        
\def\bu{\mathbf u}       

\def\Gr{Gr}     
\renewcommand\L{\mathcal{L}}

\def\tr{\mathrm{tr}\,}   

\def\pr{\mathrm{pr}}     

\def\R{\mathbb R}        
\def\C{\mathbb C}        
\def\Ad{\mathrm{Ad}}     
\def\CP{\mathbb{C}\mathrm{P}}   
\def\SU{\mathrm{SU}}     
\def\SO{\mathrm{SO}}     
\def\U{\mathrm{U}}       
\def\su{\mathfrak{su}}   
\def\SL{\mathrm{SL}}     
\renewcommand{\sl}{\mathfrak{sl}}   

\def\Re{\mathrm{Re}}  
\def\Im{\mathrm{Im}}  
\def\Imk{\mathrm{Im}\langle \cdot,\cdot \rangle}

\def\span{\mathrm{span}} 

\newtheorem{lemma}{Lemma}
\newtheorem{thm}{Theorem}
\newtheorem*{thm*}{Theorem}
\newtheorem*{notation}{Notation}
\newtheorem{corollary}{Corollary}
\newtheorem{prop}{Proposition}

\def\zbar{\overline{z}}                        
\def\delz #1{\frac{\partial}{\partial z_{#1}}} 
\def\delzbar #1{\frac{\partial}{\partial \zbar_{#1}}} 
\def\delw{\frac{\partial}{\partial w}}
\def\delwbar{\frac{\partial}{\partial \overline{w}}}

\begin{document}

\title{Compact Symmetric Spaces, Triangular Factorization, and Poisson Geometry}
\author{Arlo Caine\footnote{Research supported in part through the NSF \emph{VIGRE} Grant at the University of Arizona}
 \\Department of Mathematics \\University of Arizona \\ caine@math.arizona.edu }
\date{}
\maketitle

\begin{abstract}
Let $X$ be a simply connected compact Riemannian symmetric space, let $U$ be the universal covering group
of the identity component of the isometry group of $X$, and
let $\g$ denote the complexification of the Lie algebra of $U$, $\g=\u^\C$.  Each $\u$-compatible triangular decomposition
$\g=\n_- + \h + \n_+$ determines a Poisson Lie group structure $\pi_U$ on $U$.
The Evens-Lu construction (\cite{EL}) produces a $(U,\pi_U)$-homogeneous Poisson structure on $X$.
By choosing the basepoint in
$X$ appropriately,  $X$ is presented as $U/K$ where $K$ is the fixed point
set of an involution which stabilizes the triangular decomposition of $\g$.
With this presentation, a connection is established between
the symplectic foliation of the Evens-Lu Poisson structure
and the Birkhoff decomposition of $U/K$.  This
is done through reinterpretation of results of Pickrell (\cite{pickrell}).
Each symplectic leaf admits a natural torus action.  It is shown that the action is
Hamiltonian and the momentum map is computed using triangular
factorization.  Finally, local formulas for the
Evens-Lu Poisson structure are displayed in several examples.
\end{abstract}

\section{Introduction}
Let $(U,\pi_U)$ be a Poisson Lie group and $\d$ be the double of its Lie bialgebra.
Let $\L(\d)$ denote the variety of Lagrangian subalgebras of $\d$.
Drinfeld showed in \cite{Drinfeld} that the $U$-equivariant isomorphism classes of
$(U,\pi_U)$-homogeneous Poisson spaces with connected stability subgroups are in
one-to-one correspondence with the $U$-orbits in a certain subset of $\L(\d)$.
Evens and Lu gave a general construction in \cite{EL} which produces a
$(U,\pi_U)$-homogeneous Poisson structure on each $U$-orbit in $\L(\d)$.

Let $X$ be a connected, and simply connected, compact Riemannian symmetric space.  Let $U$ be the
universal covering group of the identity component of the isometry group of $X$
and denote by $\pi_U$ a Poisson Lie group structure on $U$.  The complexification of
$U$ will be denoted by $G$ with corresponding Lie algebra $\g=\u^\C$.
This paper concerns $(U,\pi_U)$-homogeneous Poisson structures on $X$, i.e.,
structures for which the action map $U\times X\to X$ is a Poisson map.
By selecting a basepoint, which determines a stability subgroup $K$ in $U$ and an
involution $\theta$ of $U$ fixing $K$, the symmetric space $X$ is
presented as a coset space $U/K$.
This choice determines a model point $\g_0$ in $\L(\d)$ for $X$.
The Evens-Lu construction generates a $(U,\pi_U)$-homogeneous Poisson
structure on $X$.

The difference of any two
$(U,\pi_U)$-homogeneous Poisson structures on a $U$-homogeneous space $M$ is a $U$-invariant
bivector field on $M$.  If $M$ is an irreducible compact Hermitian symmetric space,
then $M$ admits a one parameter family of $U$-invariant bivector fields.
The elements of this family are the scalar multiples of the
non-degenerate $U$-invariant Poisson
structure $\pi_{KKS}$ that a Hermitian symmetric space carries
because it is a coadjoint orbit.  The structure $\pi_{KKS}$ is the
contravariant version of the symplectic structure on coadjoint orbits
discovered by Kostant, Kirillov, and Souriau.  If $M$ is a non-Hermitian
irreducible symmetric space, then only the trivial bivector field is
$U$-invariant.  Thus, there is exactly one $(U,\pi_U)$-homogeneous
Poisson structure on an irreducible non-Hermitian symmetric space.
For standard Poisson Lie group structures $\pi_U$, the
$(U,\pi_U)$-homogeneous Poisson structures on irreducible Hermitian
symmetric spaces were classified in \cite{KKR}.

Each choice of a $\u$-compatible triangular decomposition of $\g$
determines a standard Poisson Lie group structure $\pi_U$ on $U$. In this
paper the $(U,\pi_U)$-homogeneous Poisson structure $X$ resulting
from the Evens-Lu construction is studied.
The key point of this paper is that one can choose the basepoint of
$X$ in such a way that the
triangular decomposition of $\g$ defining $\pi_U$ is stable
with respect to the involution selecting the stability subgroup $K$.
With such a presentation
of $X$ as $U/K$ there is a beautiful connection between the
symplectic foliation of the
Evens-Lu Poisson structure on $X$ and the corresponding
Birkhoff decomposition of $U/K$.
This connection enables the explicit description of the symplectic
foliation.  Each symplectic leaf has a
natural Hamiltonian torus action.  The connection with the
Birkhoff decomposition enables the computation of the momentum map.
It should also be emphasized that this presentation is very useful for
doing explicit calculations to produce examples.

This Poisson structure on compact symmetric spaces was also studied by Foth and Lu
in \cite{FL}.  They gave an alternate construction which one may interpret as follows.
It is possible to choose the basepoint of $X$, determining a different stability subgroup
$K'$ and corresponding model point $\g_0'$, in such a way that the Borel subalgebra
$\h + \n_+$ is Iwasawa relative to the noncompact real form $\g_0'$ of $\g$.
This means that the intersection
$(\h+ \n_+)\cap \g_0'$ contains $\a_0'+\n_0'$ for some Iwasawa decomposition $\g_0'=\k'+\a_0'+\n_0'$ of $\g_0'$.  In this presentation,
the push-forward of $\pi_U$ under the projection map $U\to U/K'$ gives the
Evens-Lu Poisson structure.  This method of construction has the advantage that the
natural quotient map is Poisson but the drawback that it is difficult
to explicitly calculate examples.   Also, it is not clear that there is a
torus action on the maximal symplectic leaves.

It is important to note that the stability subgroup $K'$ is not
necessarily a Poisson Lie subgroup of $(U,\pi_U)$,
even though the projection of $\pi_U$ defines a $(U,\pi_U)$-homogeneous
Poisson structure on $U/K'$.  Using results of Lu and Weintstein (\cite{LW}),
the authors Khoroshkin, Radul, and Rubtsov proved the existence in
the Hermitian symmetric case of a parabolic subgroup $P$ of $G$ such that
$U\cap P$ is a Poisson Lie subgroup of $(U,\pi_U)$.  With $X$ presented as
$U/(U\cap P)$, the natural projection of $\pi_U$ to $U/(U\cap P)$ defines
another $(U,\pi_U)$-homogeneous Poisson structure on $X$. This structure
will be denoted $\pi_{PL}$ as in \cite{KKR}.  Khoroshkin, Radul, and
Rubtsov proved that the Schouten bracket $[\pi_{PL},\pi_{KKS}]$ vanishes.
Thus,
$\pi_{PL}+\lambda \pi_{KKS}$ is the one parameter family of
$(U,\pi_U)$-homogeneous Poisson structures on a Hermitian symmetric space
$X$.

This paper is organized as follows.  In section \ref{PoissonGeometry}
it is shown that each $\u$-compatible triangular decomposition of $\g$
determines a Poisson Lie group structure $\pi_U$ on $U$. The Evens-Lu
construction is reviewed producing a $(U,\pi_U)$-homogeneous Poisson
structure on $X$. In section \ref{strata}, the Birkhoff decomposition of
$X$ is reviewed and it is proven that one can choose a basepoint in $X$
in such a way that the corresponding Birkhoff decomposition aligns
with the symplectic foliation of $X$.  By
reinterpreting the results on the Birkhoff decomposition of $X$ in
\cite{pickrell}, the leaves of the symplectic foliation are characterized
and tori acting on the symplectic leaves are determined.  In section
\ref{momentum} it is shown that the tori act in a Hamiltonian fashion and
the momentum maps are computed.  Section \ref{comments} contains a
finer discussion for the cases when $X$ is a compact Lie group or the
involution $\theta$ is an inner automorphism.  Finally, in section
\ref{examples}, local formulas for the Evens-Lu Poisson structure are
displayed.  This is done for the complex Grassmannian, complex
projective space, and the compact Lie group $\SU(2)$.

\section{Poisson Geometry \label{PoissonGeometry}}
Let $X$ be a simply connected compact symmetric space.  For simplicity, further
assume that $X$ is irreducible.  Let $U$ be the universal covering group
of the identity component of the isometry group of $X$.  The
complexification of $U$ will be denoted $G$ with corresponding
Lie algebra $\g=\u^\C$.

Each choice of a $\u$-compatible triangular decomposition of $\g$ determines a Poisson Lie
group structure on $U$ as follows.  Fix a $\u$-compatible triangular decomposition
$\g=\n_-+ \h +\n_+$.  This means that: $\h$ is a Cartan subalgebra of $\g$;
a set of positive roots for the adjoint action of $\h$ on $\g$ has been chosen;
$\n_\pm$ is the direct sum of the positive (respectively negative) root
spaces; and $-(\n_\pm)^*=\n_\mp$ where $-(\cdot)^*$ denotes the Cartan involution
selecting $\u$ in $\g$.  The sum $\g=\n_-+\h+\n_+$ is a direct sum of
vector spaces.  As every sum of vector spaces that will need
to be written down in this paper will be direct, the notation $+$ will be
used instead of the more cumbersome $\oplus$ to denote direct sum.   Set $\t=\h\cap \u$ and $\h_\R=i\t$.
Define a $\C$-linear transformation $\H\colon \g\to\g$ relative
to the given triangular decomposition by
\begin{equation}
\H(Z_-+Z_\h+Z_+)=-iZ_-+iZ_+ \label{HilbertTransform}
\end{equation}
for each $Z=Z_-+Z_\h+Z_+\in \g=\n_- + \h + \n_+$.  The real subspace $\u$ of $\g$ is
stable under $\H$.  Thus $\H$ is the complexification of a
$\R$-linear transformation $\H_\R\colon \u\to\u$ which is skew-symmetric relative to the
Killing form on $\u$.  Using the Killing form, denoted $\langle \cdot, \cdot \rangle$,
identify the dual of $\u$ with $\u$ itself and then view $\H_\R$ as an element of $\u\wedge\u$.
The bivector field
\begin{equation}\label{standardPS}
\left.\pi_U\right|_g=\ell_{g*}\H_\R - r_{g*}\H_\R
\end{equation}
defines a Poisson Lie group structure on $U$. (Here $\ell_g$ and $r_g$ denote left and right
translation by $g\in U$).  The Lie algebra structure induced by $\pi_U$ on the dual of $\u$ is
isomorphic to the real Lie algebra $\n_-+\h_\R$.  The identification is
given by the imaginary part of the Killing form.  The double of this Lie bialgebra can then be
identified with $\g$ (regarded as a real Lie algebra) via the Iwasawa decomposition
$\g=\n_-+\h_\R + \u$.
The Poisson Lie group structure $\pi_U$ corresponds to the Manin
triple $(\g,\u,\n_-+\h_\R)$.  Following terminology common in the literature,
the Poisson structure in (\ref{standardPS}) will be referred to as a \emph{standard Poisson Lie group
structure} on $U$.

More often, this Poisson structure is presented in the literature in terms of a basis.  For
each positive root $\alpha$ of the action of $\h$ on $\g$ choose a root vector $E_\alpha$
such that $\langle E_\alpha, - (E_\alpha)^*\rangle=-1$.  The $\u$-compatibility
of the triangular decomposition implies that $E_\alpha^*$ is a root vector
for $-\alpha$.  Set $E_{-\alpha}=E_\alpha^*$.  Then
$X_\alpha=E_\alpha-E_{-\alpha}$ and $Y_\alpha=i(E_\alpha+E_{-\alpha})$
are in $\u$ for each positive root and
$\H_\R=\frac{1}{2}\sum_{\alpha >0} X_\alpha\wedge Y_\alpha$.
It is shown in \cite{LW} that $\H_\R$ satisfies the Yang-Baxter equation
and thus, by a theorem of Drinfeld, the bivector field $\pi_U$ in
(\ref{standardPS}) defines a Poisson Lie group structure on the compact
group $U$.

For the remainder of this paper, fix a triangular decomposition of $\g$
and thus a standard Poisson Lie group structure $\pi_U$ on $U$.

Write $\d$ for $\g$ regarded as a real Lie algebra.  The imaginary
part of the Killing form for $\g$, denoted $\Imk$, gives a real
bilinear form on $\d$.  A Lie
subalgebra $\l$ of $\d$ is said to be \emph{Lagrangian} with respect to
$\Imk$, if $2\dim_\R \l=\dim_\R \g$, and
$\Im \langle a,b\rangle =0$ for all $a,b\in \l$.   The set of all Lie
subalgebras of $\d$ which are Lagrangian
with respect to $\Imk$ will be denoted $\L(\d)$.  This is
naturally a subvariety of the real Grassmannian  $\Gr(d,\d)$ of $d$-dimensional
subspaces of $\d$, where $2d=\dim_\R \d$.  The adjoint action of $G$ on $\g$
induces an action of $G$ (and therefore any subgroup of $G$) on $\L(\d)$.
Each $U$-orbit in $\L(\d)$ is smooth.

Using $\pi_U$, Evens and Lu (\cite{EL}) construct a smooth bivector field $\Pi$ on $\Gr(d,\d)$
with the property that the Schouten bracket $[\Pi,\Pi]$ vanishes at each point of the
subvariety $\L(\d)$.  Furthermore, they show that $\Pi$ is tangent to each $U$-orbit in $\L(\d)$, so
that each $U$-orbit is a Poisson manifold.  The construction is carried out in such a way that
each $U$-orbit in $\L(\d)$ becomes a $(U,\pi_U)$-homogenous Poisson space.

Each choice of basepoint in $X$ determines a model point in $\L(\d)$ as follows.  Let $K$ be the
stability subgroup in $U$ of a basepoint in $X$.  The subgroup $K$ is closed and is the fixed
point set of an
involution $\theta$ of $U$.  Let $\theta$ also denote the complex extension of the involution
from $U$ to all of $G$ and write $g^\theta$ for $\theta(g)$.  The Cartan involution
of $G$ fixing $U$ will be denoted $g\mapsto g^{-*}$.
Since $(\cdot)^{-1}$, $(\cdot)^{-*}$, and $\theta$ all commute, the practice
of writing $\theta$ as a superscript will not cause
confusion.    Write $G_0$ for the connected subgroup of $G$ which is fixed
by the involution $g\mapsto g^\sigma=g^{-*\theta}$.  The intersection
of $U$ and $G_0$ in $G$ is $K$ and the Lie algebra of $G_0$,
denoted $\g_0$, is a real form of $\g$ and thus a Lagrangian
subalgebra of $\d$.  The coset space $U/K$ is a finite sheeted covering
of the $U$ orbit through $\g_0$ in $\L(\d)$.  This is why $\g_0$ is
called a model point for $X\simeq U/K$ in $\L(\d)$.

The diagram of groups and Lie algebras shown in figure 1 lists
this information for reference. The upward arrows
are inclusions in both diagrams.  In the group diagram the
quotients are listed for each leg.  Each quotient is also a symmetric
space.  The quotient $G_0/K$ is a model for the non-compact symmetric
space dual to $X$ presented as $U/K$.  Also shown in the diagram
is the decomposition of the Lie algebras $\u$ and $\g_0$ into
the eigenspaces of $\theta$ as $\g_0=\k + \p$ and $\u=\k + i\p$.  At times
in this paper it will be convenient to write $X^g$ for $\Ad(g)(X)$ to
compactify notation.  Unfortunately, the adjoint action of $g$ does not
necessarily commute with the involutions which are also being written
as superscripts.  Thus one notation or the other will be used
at different points in the paper depending on the situation.

\begin{figure}[t]\label{diagrams}
\[\xymatrix{
    & \g & & & & G & \\
 \g_0=\k + \p \ar[ur] & & \u=\k + i\p \ar[ul] & & G_0 \ar[ur]^{G/G_0} & & U \ar[ul]_{G/U} \\
     & \k \ar[ul] \ar[ur] &  & & &  K \ar[ul]^{G_0/K} \ar[ur]_{U/K} &
}
\]
\caption{Algebras, groups, and quotients.}
\end{figure}
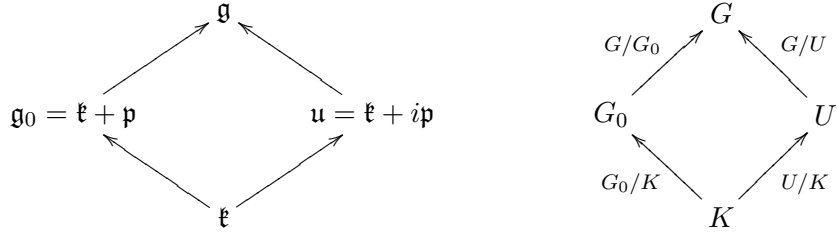

The Poisson structure $\Pi$ on the orbit $U\cdot\g_0$ may be lifted to a Poisson structure
$\pi$ on $X$ and it is this structure, making $X$ into a $(U,\pi_U)$-homogeneous Poisson space,
that it of interest in this paper.  It will be referred to as the \emph{Evens-Lu
Poisson structure} on $X$.  The Evens-Lu construction will now be reviewed
and an expression will be derived for $\pi$ amenable to the discussion
in later sections.

Let $(\cdot,\cdot)$ denote the non-degenerate $\R$-bilinear form on
$\d=(\n_- + \h_\R) + \u$ defined by
\begin{equation*}
(\xi_1+x_1,\xi_2+x_2)=\Im\langle \xi_1,x_2\rangle + \Im\langle \xi_2,x_1\rangle
\end{equation*}
for each $\xi_k \in (\n_- + \h_\R)$ and $x_k\in \u$.  Identify the dual
and double dual of $\d$ with $\d$ itself using this non-degenerate
pairing.  This
allows one to define an element of $\bigwedge^2\d$ by its action on a pair of
elements of $\d$.  Evens and Lu define
$R\in \bigwedge^2\d$ by
\begin{equation}\label{Rmatrix}
R(\xi_1+x_1,\xi_2+x_2)=\Im\langle \xi_2,x_1\rangle -\Im\langle \xi_1,x_2\rangle
\end{equation}
and use it to generate a bivector field on the Grassmannian of $d$-dimensional subspaces
of $\d$.  The adjoint action of $G$ on $\g$ induces a $\R$-Lie algebra anti-homomorphism
\begin{equation*}
\kappa \colon \textstyle{ \d \to \Gamma(T\Gr(d,\d))}
\end{equation*}
whose multi-linear extension $\bigwedge^\cdot \d\to \Gamma(\bigwedge^\cdot T\Gr(d,\d))$
will also be denoted by $\kappa$.  The bivector field $\Pi$ is defined by
$\Pi=\frac{1}{2}\,\kappa(R)$.

To do calculations, identify $T(U/K)$ with $U\times_K i\p$ using right translation.
Further, identify $i\p$ with $(i\p)^*$ using the Killing form so that
$T^*(U/K)$ may also be represented by $U\times_K i\p$.  In this setting,
the action of $\pi$ on a pair of cotangent vectors
represented by classes $[u,X]$ and $[u,Y]$ in $U\times_K i\p$ may be
computed by
\begin{equation*}
\pi([u,X],[u,Y])=\langle \Omega_u X, Y\rangle
\end{equation*}
where $\Omega_u\colon i\p \to i\p$ is a skew-symmetric $\R$-linear transformation
which is $K$-equivariant in its dependence on $u$.

\begin{thm}
The Evens-Lu Poisson bivector $\pi$ can be expressed as
\begin{equation*}
\pi([u,X],[u,Y])=\langle \Omega_u X, Y\rangle
\end{equation*}
where $X,Y\in i\p$ and $\Omega_u\colon i\p \to i\p$ is given by
\begin{equation}\label{formula_for_pi}
\Omega_u(X)=\{\Ad(u^{-1})\circ \H\circ \Ad(u) (X)\}_{i\p}
\end{equation}
and $\{\cdot \}_{i\p}$ denotes the projection to $i\p$ along the decomposition
$\g=i\u+\k+i\p$.
\end{thm}
\begin{proof}
Let $\pr_\u\colon \g\to \u$ denote the projection onto $\u$ along the Iwasawa
decomposition $\g=\n_-+\h_\R+\u$.  In terms of the triangular decomposition
$Z=Z_-+Z_\h+Z_+\in \g$, this projection is given by the formula
\begin{equation*}
\pr_\u (Z) = -(Z_+)^*+Z_\t+Z_+.
\end{equation*}
The composition $\u\xrightarrow{i} i\u \xrightarrow{\pr_\u} \u$ agrees with
the transformation $\H$ when restricted to $\u$ (recall the definition of
$\H$ in (\ref{HilbertTransform})).  Indeed, each $Z\in \u$
satisfies $Z_-=-(Z_+)^*$ and $Z_\h=Z_\t$; hence $(iZ)_\t=0$,
and $-((iZ)_+)^*=i(Z_+)^*=iZ_-$. It follows that
\begin{equation}\label{i_proj_agrees_with_H}
\pr_\u (iZ)=-iZ_-+iZ_+ = \H (Z).
\end{equation}

In equation (\ref{formula_for_pi}) the equivalence class $[u,X]$ represents
a linear functional
on the tangent space to $U/K$ at $uK$.  Using the covering map to identify
this space with the tangent space to the $U$-orbit through $\g_0$ at $u\cdot \g_0$,
each such tangent vector corresponds to an element $\chi\in \u/\k^u$ via the
map $\chi\mapsto \left.\frac{d}{dt}\right|_{t=0} e^{t\chi} u\cdot \g_0$.  The
corresponding class in $U\times_K i\p$ is $[u,\chi^{u^{-1}}]$.  Moreover,
the action of $[u,X]$ as a linear functional on $[u,\chi^{u^{-1}}]$ is given by
\begin{equation*}
[u,X]\left([u,\chi^{u^{-1}}]\right)=\langle X,\chi^{u^{-1}}\rangle =\Im \langle iX^u,\chi\rangle .
\end{equation*}
By the Evens-Lu construction
\begin{equation*}
\pi ([u,X],[u,Y]) =\frac{1}{2}\left.\kappa (R)\right|_{uK} ([u,X], [u,Y]) = \frac{1}{2} R(iX^u,iY^u) .
\end{equation*}
Combining the definition of $R$ in (\ref{Rmatrix}), the result in (\ref{i_proj_agrees_with_H}),
and skew-symmetry of $\H$, it follows that
\begin{eqnarray*}
R(iX^u,iY^u) & = & \Im \langle \pr_\u (iX^u),iY^u\rangle - \Im \langle iX^u, \pr_\u (iY^u)\rangle \\
& = & \langle \H (X^u), Y^u\rangle -\langle X^u , \H (Y^u)\rangle \\
& = & 2\langle (\H (X^u))^{u^{-1}}, Y\rangle .
\end{eqnarray*}
This completes the proof of the theorem.
\end{proof}
From this formula, the connection with triangular factorization is evident and
a group theoretic interpretation of the symplectic foliation can be given.
At the group level write the corresponding Iwasawa decomposition for $G$
as
\begin{equation}
\begin{array}{rcl}G & \simeq & N^- \times A \times U .\\
                  g & \mapsto &  (\mathbf{l}(g),\mathbf{a}(g),\bu (g))
\end{array}
\end{equation}
where $A=\exp(\h_\R)$. There is a natural action of
$G$ (and therefore any subgroup of $G$) on $U$ coming from the identification of $U$
with the right coset space $N^-A\backslash G$.
\begin{equation*}
\begin{array}{rcl} U\times G & \to & U \\
                   u \cdot g &\mapsto & \bu(ug)
\end{array}
\end{equation*}
As observed in \cite{LW}, the symplectic leaves of the Poisson Lie group
structure $\pi_U$ on $U$ are precisely the $N^-A$ orbits in $U$.
\begin{prop}\label{me1}
For the Evens-Lu Poisson structure $\pi$ the leaves of the
symplectic foliation are the projections of the $G_0$-orbits
in $U$ to $U/K$.
\end{prop}
\begin{proof} The symplectic foliation of $U/K$ is generated by the distribution in $T(U/K)$
which is image of the natural map $\pi^\natural\colon T^*(U/K)\to T(U/K)$ given by
$\pi^\natural (\alpha)=\pi(\alpha , \cdot)$.  In terms of the identifications
$T^*(U/K)\simeq U\times_K i\p\simeq T(U/K)$,
\begin{equation*}
\pi^\natural ([u,X])=[u,\{(\H(X^u))^{u^{-1}}\}_{i\p}]=[u,\{ (\pr_\u (iX^u))^{u^{-1}}\}_{i\p}].
\end{equation*}
Fix $u\in U$. The map $G_0\to U$ given by $g_0\mapsto \bu(ug_0)$ is equivariant
for the right $K$-action on both $G_0$ and $U$ and thus descends to a
map $G_0/K\to U/K$.
Now, suppose that $iX \in \p$, so that $X\in i\p$.  Consider the curve
\begin{equation*}
u\cdot e^{tiX} = \bu(ue^{tiX})=\bu(e^{tiX^u}u).
\end{equation*}
Differentiating at $t=0$, one obtains the tangent vector
$(r_u)_*(\pr_\u(iX^u))$.
It follows that the distribution tangent to the projection of the
$G_0$-orbit at $uK$
and is spanned by vectors of the form
\begin{equation}
\kappa (\pr_\u (iX^u))|_{uK}=[u,\{(\pr_\u (iX^u))^{u^{-1}}\}_{i\p}]
\end{equation}
which is exactly the image of $\pi^\natural$.
\end{proof}
This proposition has been established in several contexts for this Poisson
structure.  Using their construction, Foth and Lu give an alternate proof
(cf. Proposition 1.1 in \cite{FL}).

\section{Symplectic Leaves and Triangular Factorization\label{strata}}

Let $G$ be a simply connected complex semi-simple Lie group with Lie
algebra $\g$.  Choose a Cartan subalgebra $\h$, and a set of positive roots
for the adjoint action of $\h$ on $\g$.  Let $\n_\pm$ denote the sum of the
positive (resp. negative) root spaces.  This data gives a triangular
decomposition of $\g$.
\begin{equation*}
\g = \n_- +\h +\n_+
\end{equation*}
Set $H=\exp(\h)$, and $N^\pm=\exp(\n_\pm)$.  Corresponding
to this decomposition is a Birkhoff decomposition (a.k.a. triangular
decomposition a.k.a LDU decomposition) of the group
\begin{equation*}
G=\coprod_{w\in W} \Sigma_\w^G\text{, where }\Sigma_w^G= N^-wHN^+
\end{equation*}
where $w$ is an element of the Weyl group $W=N_G(H)/H$.  Each component
$\Sigma_w^G$ is a manifold diffeomorphic to $N^-\cap wN^-w^{-1}\times H \times N^+$.
The codimension of $\Sigma_w^G$ increases with the length of $w$ in the Weyl
group and $\Sigma_1^G$ is a Zariski open subset of $G$.  Each element in
$\Sigma_w^G$ can be factored as a product of an element of $N^-$, an
element of $w\subset N_G(H)/H$, and element of $H$, and an element of $N^+$.
For $\SL(n,\C)$ in an appropriate representation this would correspond
to the factorization of an $n\times n$ complex matrix of determinant one into a
product of a lower triangular unipotent matrix, a unitary permutation matrix,
a diagonal matrix of determinant one, and an upper triangular unipotent matrix.
Each element of $\Sigma_1^G$ admits a unique factorization of this form.
If $w$ is a non-trivial element of the Weyl group then the elements in
$\Sigma_w^G$ admit several such factorizations.  Further conditions
are required to guarantee uniqueness in those cases.

This section concerns a generalization of this decomposition for symmetric
spaces and, in particular, its relationship to the Poisson geometry of $(X,\pi)$.
Recall the setup considered in this paper: $X$ is a compact, connected, and
simply connected Riemannian symmetric space; $U$ is the universal covering
group of the identity component of the isometry group; $G$ is the complexification
of $U$.  The Poisson Lie group structure $\pi_U$ on $U$ was defined by
a triangular decomposition of $\g$ which had the additional property of
being $\u$-compatible, i.e. $-(\n_\pm)^*=\n_\mp$.

For each choice of basepoint in $X$, Cartan defined an embedding of $X$ into
$U\subset G$.
\begin{equation}\label{cartan_embedding}
\begin{array}{rccccc}
\phi \colon & U/K &   \to    &      U       & \hookrightarrow & G \\
            & uK  & \mapsto  & uu^{-\theta} &                 &   \\
\end{array}
\end{equation}
In (\ref{cartan_embedding}), $\theta$ is the involution fixing the stability
subgroup of the basepoint in $X$. If $U$ is viewed as a Riemannian
manifold with the metric induced from the Killing form, then this map
gives a totally geodesic embedding of $X$ into $U\subset G$. The
intersection of this image with the decomposition of $G$ induces
a decomposition of $X$ which will be called a \emph{Birkhoff decomposition
of $X$}.
Such a decomposition depends on the choice of basepoint in
$X$ and the triangular decomposition of $\g$.  In this paper, the triangular
decomposition of $\g$ is regarded as fixed, having defined $\pi_U$.

\begin{lemma}\label{basepoint_choice}
There exists a basepoint in $X$ such that the given triangular
decomposition of $\g$ is stable with respect to the involution selecting
the stability subgroup of $x$ in $U$.
\end{lemma}
\begin{proof}
Fix a point $x\in X$.  This determines the data: $K$, the stability
subgroup of $x$; $\theta$, the involution of $U$ fixing $K$; $G_0$, a non-compact real
form $G$; $\u=\k+i\p$, a decomposition of $\u$ into the eigenspaces of $\theta$.
The Lie algebra of $G_0$ is $\g_0=\k+\p$.  First, it will be shown that a
triangular decomposition of $\g$, which is both $\u$-compatible
and $\theta$-stable, exists.

Let $\t_0$ be a Cartan subalgebra of $\k$, and $\h_0$ denote the centralizer
of $\t_0$ in $\g_0$.  Theorem 6.60 of \cite{Knapp} shows that $\h_0=\t_0+\a_0$
is a Cartan subalgebra of $\g_0$, where $\a_0\subset \p$.  Thus, the subalgebras
\begin{equation*}
\h_0=\t_0+\a_0\text{, }\t=\t_0+i\a_0\text{, and }\h=\t^\C
\end{equation*}
are Cartan subalgebras of $\g_0$, $\u$, and $\g$ respectively.  Selecting
a Weyl chamber in $\h_\R=i\t$, chooses a set of positive roots for the
adjoint action of $\h$ on $\g$, and thus determines a triangular decomposition
for $\g$.  For each choice the resulting decomposition will be $\u$-compatible.
The Cartan subalgebra $\h$ is $\theta$-stable by construction.
However, $\theta$ will stabilize the positive root spaces only if the chosen
Weyl chamber contains an element of $i\t_0$.  The fact that $i\t_0$ contains
a regular element is equivalent to the fact that $\h_0$ is a Cartan
subalgebra of $\g_0$.  Thus, a $\u$-compatible $\theta$-stable triangular
decomposition of $\g$ exists.

To finish the proof, note that all $\u$-compatible triangular decompositions of $\g$
are $U$-conjugate.  Therefore, one may conjugate the constructed
triangular decomposition to the given $\u$-compatible
triangular decomposition of $\g$ by an element $u\in U$.  Conjugating the
stability subgroup $K$ by the same element selects a new basepoint in $X$
with corresponding involution $\theta'=\Ad(u)\circ \theta \circ \Ad(u^{-1})$.  The
given $\u$-compatible decomposition $\g=\n_-+\h+\n_+$ is stable with
respect to the involution $\theta'$.
\end{proof}

For the remainder of this section, fix a presentation of the symmetric space
$X$ as $U/K$ such that the triangular decomposition of $\g$ defining
$\pi_U$ is $\theta$-stable.  Under the assumption of $\theta$-stability,
Pickrell was able to characterize the connected
components of the resulting Birkhoff decomposition of $X$
(cf. \cite{pickrell}).  It is
through reinterpretation of these results that an explicit
description of the geometry of the Evens-Lu Poisson structure on
$X$ can be given.  For the convenience of the reader, the same notation
as that in \cite{pickrell} has been adopted in this paper.

The first observation is that $\phi(U/K)$ does not necessarily intersect
each component of the Birkhoff decomposition of $G$.  In fact,
$\phi(U/K)\cap \Sigma_w^G$ is non-empty if and only if $\phi(U/K)\cap w$ is
non-empty where $w\subset N_U(T)$ represents the Weyl group element
$w\in W=N_G(H)\simeq N_U(T)$ (cf. Theorem 2 (a) combined with Theorem 1 (a) in
\cite{pickrell}).

One should think of a Birkhoff decomposition of $X$ as consisting of a
number of layers indexed by the possible elements $w$.  In this paper,
\emph{the layer corresponding to $w$}, will refer to the set
$\phi(U/K)\cap \Sigma_w^G$, or
its equivalent in $U/K$ or $X$ when the context is clear.  The reader
should understand that when this terminology is used they are
to implicitly assume that $w$ is such that $\phi(U/K)\cap w$ is non-empty.

Each layer consists of a number of connected
components.  For a given $w\in N(T)/T$ the connected components
of the layer corresponding to $w$ are indexed by the set
\begin{equation}
\{\w\in \phi(U/K)\cap w\}/T
\end{equation}
where $T$ acts on the right by $\w\cdot t=t^{-1}\w t^\theta$ (cf. Theorem 2 (c) in \cite{pickrell}).
This characterization uses $\theta$-stability.  The elements $\w\in \phi(U/K)\cap T$
are the images under the Cartan embedding of the preferred basepoints
in $X$ in the sense of lemma \ref{basepoint_choice}.

\begin{notation} \emph{
Given an element $\w$ in the layer corresponding to $w$, write
$\Sigma_\w^{\phi(U/K)}$ for the connected component of $\phi(U/K)\cap \Sigma_w^G$
containing $\w$.}
\end{notation}

In proposition \ref{me1} it was shown that the symplectic leaves
the Evens-Lu Poisson structure on $X\simeq U/K$ are the projections modulo $K$ of the
$G_0$-orbits in $U$.  The right action of $G_0$ on $U$ was that induced
by the identification $U\simeq N^-A\backslash G$ coming from the Iwasawa
decomposition $G\simeq N^-AU$.  By combining this with an action of
the torus $T=\exp(\h\cap \u)$ as
\begin{equation*}
\begin{array}{ccc}
 U\times (T\times G_0) &   \to   &        U        \\
    u\cdot (t,g_0)     & \mapsto & t^{-1}\bu(ug_0) \\
\end{array}
\end{equation*}
Pickrell was able to characterize each component $\Sigma_\w^\phi(U/K)$.
The following proposition is Theorem 4 (a) in \cite{pickrell}.
\begin{prop}\label{pickrell2}
Consider the layer of the Birkhoff decomposition
of $X$ corresponding to $w\in N_U(T)/T$.  Let $\w\in \phi(U/K)\cap w$ and fix a choice
of $\w_1\in U$ such that $\phi(\w_1K)=\w$. The map
\begin{equation*}
\begin{array}{ccc}
 T\times G_0 &   \to   &  \phi(U/K)\cap \Sigma_\w^{\phi(U/K)}  \\
   (t,g_0)   & \mapsto & \phi(t^{-1}\bu (\w_1 \,g_0)) \\
\end{array}
\end{equation*}
is surjective and induces a diffeomorphism
\begin{equation*}
T\times_{\exp(\ker\{\Ad(w)\theta|_\t-1\})} R\backslash G_0/K \to \Sigma_\w^{\phi(U/K)}
\end{equation*}
where $R=(N^- A)^{\w_1^{-1}}\cap G_0$ is a contractible subgroup of $G_0$
and $\lambda \in \exp(\ker\{\Ad(w)\theta|_\t-1\})$
is identified with the pair $(\lambda, \lambda^{\w_1^{-1}})$.
\end{prop}
The proof is fairly involved and will not be reproduced here.  The key
point is that the theorem can be reinterpreted as characterizing the
symplectic leaves of the Evens-Lu Poisson
structure.  The projection modulo $K$ of the $G_0$-orbit through $\w_1$
maps to a sub-manifold of $\Sigma_\w^{\phi(U/K)}$ passing through $\w$.
This shows that with the appropriate choice of basepoint in $X$, the
symplectic foliation of the Evens-Lu Poisson structure aligns with
the Birkhoff decomposition of $X$.
\begin{corollary}\label{corollary_on_symplectic_leaves}
The sorting of the symplectic leaves in $(X,\pi)$ by dimension corresponds to the
Birkhoff decomposition of $X$.  In particular, each connected component of
the layer corresponding to $w$ is foliated by contractible symplectic leaves
each diffeomorphic to the double coset space $R\backslash G_0/K$.
\end{corollary}

More can be seen from proposition \ref{pickrell2}.  Each leaf admits a
natural torus action.  The acting torus is a subgroup of the torus $T$
determined by the layer of the Birkhoff decomposition in which the leaf
is contained.

\begin{corollary}
Each symplectic leaf foliating the layer corresponding to $w$ is acted on
by the torus
\begin{equation*}
T_w = \exp(\ker \{\Ad(w)\theta|_\t -1\})\subset T.
\end{equation*}
Moreover, the coset $\w_1K$ represents the unique fixed point for the
$T_w$-action in the symplectic leaf through $\w_1K$.
\end{corollary}

In the next section, it will be shown that the action of $T_w$ on each
symplectic leaf foliating the layer indexed by $w$ is Hamiltonian.  The
following proposition and theorem will be needed in order to compute the
momentum map. First, notice that an element $g$ in $\phi(U/K)\cap G$ first
satisfies the equation $g^{-1}=g^\theta$ and the further condition that $g\in U$
so that $g^{-1}=g^*$.  In fact $\phi(U/K)$ can be realized as the identity
component of the set $\{g^*=g^\theta\}\cap U$ (cf. Theorem 1 (a) in \cite{pickrell}).  The following proposition
is Theorem 2 (d) in \cite{pickrell}.

\begin{prop}\label{pickrell_factor_theorem}
Fix $w\in W$.  Suppose $\w\in w\subset N_G(H)$  with $\w^*=\w^\theta$.
For brevity, write $D_\w$ for the set
\begin{equation*}
\{(h,L)\in H\times (N^-\cap (N^+)^w)\colon \theta(h^{w^{-1}})=h^*,
(\theta(L^{-*}))^{\w h}=L^{-1}\}.
\end{equation*}
The map
\begin{equation*}
\begin{array}{ccc}
( N^-\cap (N^-)^w)\times D_\w &   \to   & \{g^*=g^\theta\}\cap \Sigma_w^G \\
        (l,(h,L))        & \mapsto & lL^{-1}\w h (lL^{-1})^{*\theta} \\
\end{array}
\end{equation*}
is a diffeomorphism onto the connected component of $\{g^*=g^\theta\}\cap \Sigma_w^G$
containing $\w$.
\end{prop}

Applied to $\phi(U/K)\cap \Sigma_w^G$, this proposition provides
further conditions to guarantee the uniqueness of the triangular factorization
of elements in the layers indexed by non-trivial elements $w$.  The following
corollary reinterprets this result.

\begin{corollary}\label{me2}
Let $\w\in \phi(U/K)\cap w$ be an element in the layer corresponding to $w$.
Let $S_\w$ denote the symplectic leaf of $(X,\pi)$ for which
$\phi(S_\w)$ passes through $\w$ in $\phi(U/K)\cap w$.
Each element in $\phi(S_\w)$ can be factored as $\ell \w h \ell^{*\theta}$ where
$\ell \in N^-$, and $h\in \exp(\ker\{Ad(w)\sigma|_\h +1\})$.  Furthermore, the
magnitude $|h|=\sqrt{h^*h}$ is a smooth function on the leaf.
\end{corollary}
\begin{proof}
Existence of the factorization is guaranteed by proposition
\ref{pickrell_factor_theorem} with the observation that $h$ satisfies
the condition $\theta(h^{w^{-1}})=h^*$ if and only if
$h\in \exp(\ker\{\Ad(w)\sigma|_\h +1\})$.  Note that proposition
\ref{pickrell_factor_theorem} also guarantees that $h$ is a smooth
function on the leaf and thus $|h|=\sqrt{h^*h}\in \exp(\h_\R)$ is
smooth as well.
\end{proof}

\section{The Momentum Map\label{momentum}}

Given the amount of notation introduced and the number of objects and
parameters involved it seems apropos to include a brief summary of the
setup and results of the paper thus far.  With a compact, connected,
and simply connected symmetric space $X$ comes the data: $U$, the
universal cover of the identity component of the isometry group of $X$;
$G$, the complexification of $U$; the Cartan involution $g\mapsto g^{-*}$
selecting $U$ as a compact real form of $G$.  In section \ref{PoissonGeometry}
a $\u$-compatible (i.e., $-(\n_\pm)^*=\n_\mp$) triangular decomposition
$\g=\n_-+\h+\n_+$ was fixed, which determined a standard Poisson Lie group
structure $\pi_U$ on $U$.  The Evens-Lu construction produced a $(U,\pi_U)$-homogeneous Poisson
structure $\pi$ on $X$.  In section \ref{strata} a basepoint was chosen in
$X$ whose corresponding involution $\theta$
stabilized the triangular decomposition $\g=\n_-+\h+\n_+$ ( i.e.,
$\theta(\n_\pm)=\n_\pm$, and $\theta(\h)=\h$).  This choice determined:
a presentation of $X$ as a coset
space $U/K$; an embedding $\phi\colon U/K\to U$ of $X$ into $U\subset G$.
The intersection of the image of this map with the Birkhoff decomposition
of $G$ corresponding to the fixed triangular decomposition of $\g$ induced
a Birkhoff decomposition of $X$.

In section \ref{strata} it was shown that the layers of the Birkhoff
decomposition of $X$ are foliated by symplectic leaves, each of which
is contractible.  There is a bijection between $\cup_{w\in W} (\phi(U/K)\cap w)$
and the symplectic leaves in $X$ assigning a leaf to the element its
$\phi$-image passes through.
Given $\w\in\phi(U/K)\cap \Sigma_w^G$ the corresponding leaf was denoted $S_\w$.
The torus $T=\exp(\h\cap \u)$ acts from the right on $U/K$ via
$uK\cdot t = \t^{-1}uK$.  The sub-torus $T_w=\exp(\ker\{\Ad(w)\theta|_\t-1\})$
preserves the symplectic leaf $S_\w$ with a unique fixed point.
Each element of $\phi(S_\w)$ can be factored as $\ell \w h \ell^{*\theta}$
where $\ell\in N^-$ and $h$ is uniquely determined.  The magnitude of
$h$ is a smooth function on $S_\w$.

\begin{thm}\label{me3}
The action of the torus $T_w$ on the symplectic leaf $S_\w$
is Hamiltonian with momentum map
\begin{equation*}
\begin{array}{cccc}
 \mu\colon & S_\w &   \to   &                       \t_w^*                        \\
           &  uK  & \mapsto & \langle \frac{1}{2}i\theta(\log | h |), \cdot \rangle \\
\end{array}
\end{equation*}
where $h$ is the diagonal part of $\phi(uK)=uu^{-\theta}\in \phi(S_\w)$.
\end{thm}
\begin{proof} Let $X\in \t_w$.  Practically, this means that $X$ satisfies
the equation $\w X^\theta \w^{-1} = X$ for each $\w\in w$.
The vector field $\Tilde{X}$ induced by the action $T_w$ on $U/K$ is
represented by the
class $[u,\{-X^{u^{-1}}\}_{i\p}]$ in $U\times_K i\p$. It needs to be
shown that $\mu_X\colon S_\w\to \R$ given by
$\mu_X(uK)=\langle \frac{1}{2}i\theta(\log |d_\phi |), X \rangle$ is a
Hamiltonian function for the vector field $\Tilde{X}$. This will follow
if $\pi^\natural (d\mu_X)=\Tilde{X}$.

Let $[u,Y]\in U\times_K i\p$ be a tangent vector to $U/K$ at $uK$.
This tangent vector is represented by the curve $ue^{tY}K$ which
passes through $uK$ at
$t=0$.  Decompose $\phi(uK)=uu^{-\theta}$ as
$\ell \w h \ell^{*\theta}$ using corollary \ref{me2}.  Then
\begin{eqnarray*}
ue^{tY}(ue^{tY})^{-\theta} & = & e^{2tY^u}uu^{-\theta} \\
& = & e^{2tY^u}\ell \w h \ell^{*\theta} \\
& = & \ell \w e^{2t (\ell\w)^{-1}Y^{u}\ell\w} h \ell^{*\theta}
\end{eqnarray*}
and thus $\left.\frac{d}{dt}\right|_{t=0}\log |d_\phi(ue^{tY}K)|=2\,
\pr_{\h_\R} ((\ell\w)^{-1}Y^{u}\ell \w)$.  With this calculation it follows
that
\begin{eqnarray*}
d\mu_X ([u,Y]) & = & \langle i\theta(\pr_{\h_\R} ((\ell\w)^{-1}Y^{u}\ell\w)), X\rangle \\
& = & \langle Y, \Ad(u^{-1})((i\w X^\theta \w^{-1})^{\ell}) \rangle \\
& = & \langle Y, \Ad(u^{-1})(iX^{\ell}) \rangle .
\end{eqnarray*}
Hence, $d\mu_X$ corresponds to the class
$[u, \{\Ad(u^{-1})(iX^{\ell})\}_{i\p}]$.
\medskip

\noindent {\bf Assertion 1:} $\{\Ad(u^{-1})(iX^{\ell})\}_{i\p} =
\frac{1}{2}\{\Ad(u^{-1})(iX^\ell)\}_{\u}$ where $\{\cdot \}_\u$ denotes
the projection to $\u$ along the decomposition $\g=\u+i\u$.
\medskip

Given $Z\in \g$, the projection to $i\p$ is given by
\begin{equation}
\{Z\}_{i\p}=\frac{1}{4}\left(Z+Z^{*\theta}-(Z+Z^{*\theta})^* \right)=\frac{1}{2}\{Z+Z^{*\theta}\}_\u. \label{proj_to_ip}
\end{equation}
The equation $uu^{-\theta}=\ell\w h\ell^{*\theta}$ implies
that $u^{-\theta}\ell^\sigma=u^{-1}\ell \w h$.  This equality gives
\begin{eqnarray}
(\Ad(u^{-1})(iX^{\ell}))^{*\theta} & = & u^{-\theta}\ell^\sigma iX^\theta \ell^{-\sigma} u^\theta  \label{miracle1}\\
& = & u^{-1}\ell \w iX^\theta \w^{-1} \ell^{-1} u \label{miracle2}\\
& = & \Ad(u^{-1})(iX^\ell) \label{crux}
\end{eqnarray}
where (\ref{miracle2}) follows from (\ref{miracle1}) after substituting for
$u^{-\theta}\ell^\sigma$ and noting the fact that $h \in H$ trivially
on $iX^\theta$.   The assertion follows.

In terms of the identifications of $T(U/K)$ and $T^*(U/K)$ with $U\times_K i\p$
the map $\pi^\natural$ is given by
\begin{equation} \label{pi_natural}
[u,Y]\mapsto [u,\{(\H(Y^u))^{u^{-1}}\}_{i\p}]
\end{equation}
where $\H(Y)=-i(Y)_-+i(Y)_+$ is the transformation from (\ref{HilbertTransform}).
Using assertion 1 and the fact that $\Ad(u)$ commutes with the projection
to $\u$ one can see that
$\pi^\natural (d\mu_X) = [u,\{\Ad(u^{-1})\circ \H(\{iX^\ell\}_\u)\}_{i\p}]$.

\medskip

\noindent {\bf Assertion 2:} $\H(\{iX^\ell\}_\u)= \frac{1}{2}\{X^\ell\}_{\u}-X.$
\medskip

The key point is that $iX^\ell=(iX^\ell)_-+iX$
since $\ell \in N^-$ and $X\in H$.  Hence,
\begin{equation*}
2 \{iX^\ell\}_\u = \left((iX^\ell)_- -((iX^\ell)_-)^*\right)
= \left((\ell iX\ell^{-1})_- -(\ell^{-*} iX \ell^{*})_+\right).
\end{equation*}
Applying $\H$ to this expression gives
\begin{eqnarray*}
2 \H(\{iX^\ell\}_\u) & = & (\ell X\ell^{-1})_- +(\ell^{-*} X \ell^{*})_+ \\
& = & \ell X \ell^{-1} - X + \ell^{-*}X\ell^* - X \\
& = & \{X^\ell\}_{\u} - 2X .
\end{eqnarray*}
This proves assertion 2.

To finish the calculation of $\pi^\natural (d\mu_X)$, observe that
$\Ad(u^{-1})$ preserves $\u$, and use assertion 2 to conclude that
the right hand side of (\ref{pi_natural}) is equal to
\begin{equation*}
[u, \frac{1}{2}\{\Ad(u^{-1}(X^\ell)\}_{i\p}-\{X^{u^{-1}}\}_{i\p} ].
\end{equation*}
The same calculation in (\ref{miracle1}) through (\ref{crux}) shows that
$(\Ad(u^{-1})(X^\ell))^{*\theta}=-\Ad(u^{-1})(X^\ell)$.  Combining this observation
with the formula in (\ref{proj_to_ip}) yields that the $i\p$ part of
$\Ad(u^{-1})(X^\ell)$ is zero.  Thus
\begin{equation*}
\pi^\natural (d\mu_X)=[u,-\{X^{u^{-1}}\}_{i\p}],
\end{equation*}
which completes the proof of the theorem.
\end{proof}

\section{Comments\label{comments}}

This section addresses several special cases and sets the stage for
the explicit examples in section \ref{examples}.

Let $K$ be a connected and simply connected compact Lie group.  Fix a
triangular decomposition
\begin{equation}\label{newtridecomp}
\k^\C=\tilde{\n}_-+\tilde{\h}+\tilde{\n}_+ .
\end{equation}
Write $\H_\k$ for the linear transformation in (\ref{HilbertTransform})
relative to the decomposition (\ref{newtridecomp}).  Trivialize the tangent bundle to $K$
using right translation.  This identifies $T K$ with $K\times \k$.
Further identify $\k$ will its dual using the Killing form so that $T^*K$
is also identified with $K\times \k$.  Using the definition of $\pi_K$ from
(\ref{standardPS}), a short calculation yields the following formula.
\begin{prop}
The Lu-Weinstein Poisson Lie group structure $\pi_K$ can be expressed by
\begin{equation}
\pi_K ((k,P),(k,Q))=\langle (\Ad(k)\circ \H_\k\circ \Ad(k^{-1})-\H_\k)(P),Q\rangle .
\end{equation}
for each $(k,P),(k,Q)\in K\times \k$.
\end{prop}
With this formula, one can see that the maximal torus $\t=\tilde{\h}\cap \k$
is a Poisson Lie subgroup of $K$, as $\pi_K$ vanishes identically there.
Thus, the push-forward of $\pi_K$ under the natural projection map
$K\to K/T$ defines a $(K,\pi_K)$-homogeneous Poisson structure on the
flag manifold $K/T$.  The symplectic leaves of this induced structure are
precisely the corresponding Bruhat cells in $K/T$ because the symplectic
leaves of $\pi_K$ foliate the components of the Bruhat decomposition of
$K$.

When equipped with the invariant metric induced by the Killing form, $K$
is also a Riemannian symmetric space.  The isometries of $K$ are given by either
left or right translation by elements of $K$.  In this case,
$U=K\times K$, $G=K^\C\times K^\C$, and the Cartan involution selecting
$U$ in $G$ is $(g_1,g_2)\mapsto (g_1^{-*},g_2^{-*})$ where $g\mapsto g^{-*}$
is the Cartan involution selecting $K$ inside of $K^\C$. The left action of
$u=(k_1,k_2)$ on $k\in K$ is given by $(k_1,k_2)\cdot k=k_1kk_2^{-1}$.
The decomposition
\begin{equation}\label{groupcasedecomposition}
\g= \underbrace{(\tilde{\n}_-\times \tilde{\n}_-)}_{\n_-}+
\underbrace{(\tilde{\h}\times \tilde{\h})}_{\h}+
\underbrace{(\tilde{\n}_+\times\tilde{\n}_+)}_{\n_+}
\end{equation}
built using (\ref{newtridecomp}) is a $\u$-compatible triangular
decomposition of $\g$.

Using the identity in $K$ as a basepoint gives a presentation of $K$ as a the coset space
$U/\Delta$ where the stability subgroup $\Delta=\{(k,k)\colon k\in K\}$
is the diagonal image of $K$ in $U$.  The involution $\theta$ fixing
$\Delta$ in $U$ is the automorphism which interchanges the two factors of
$U$.  The triangular decomposition in (\ref{groupcasedecomposition}) is
stable with respect to this outer automorphism.

The image of the Cartan embedding of $U/\Delta$ into $U$ is the anti-diagonal
image of $K$ in $U$, namely $\phi(U/\Delta)=\{(k,k^{-1})\colon k\in K\}\subset U$.
The corresponding Birkhoff decomposition of $K$ as a symmetric space is
identical to the decomposition of $K$ induced by the Birkhoff decomposition
of $K^\C$ with respect to $\k^\C=\tilde{\n}_-+\tilde{\h}+\tilde{\n}_+$.

\begin{thm}
The Evens-Lu Poisson structure $\pi$ on $K$ can be expressed as
\begin{equation}
\pi ((k,P),(k,Q))=\langle (\H_\k + \Ad(k)\circ \H_\k \circ \Ad(k^{-1}))(P),Q\rangle .
\end{equation}
\end{thm}
\begin{proof}
The isomorphism $\psi\colon U/\Delta \to K$ given by $(k_1,k_2)\Delta \mapsto k_1k_2^{-1}$
identifies $U/\Delta$ with $K$.  Use right translation to identify the
tangent bundle to $U/\Delta$ with $U\times_\Delta i\p$ and $TK$ with
$K\times \k$.  Given $(k_1,k_2)\in U$ write $k=\psi(k_1,k_2)=k_1k_2^{-1}$.
In this setting $i\p=\{(X,-X)\colon X\in \k\}$.  A curve
representing $(X,-X)$ in $i\p$ passing through $(k_1,k_2)\Delta$
is given by $(k_1e^{tX},k_2e^{-tX})\Delta$ and its image under $\psi$
is $k_1e^{tX}e^{tX}k_2^{-1}=e^{2tX^{k_1}}k$.  Thus,
$\psi_*[(k_1,k_2),(X,-X)]$ is represented
by $(k,2X^{k_1})\in K\times\k$.  Using the Killing form to
identify the dual of $\k$ with $\k$, the cotangent bundle can also be
represented by $K\times \k$ using right translation.  Then, the pull-back of
a cotangent vector at $k$ represented by $(k,P)$ is given by the class
$[(k_1,k_2),(P^{k_1^{-1}},- P^{k_1^{-1}})]$.  Let $\pi_{EL}$ denote the
Evens-Lu Poisson structure $U/\Delta$.  From theorem \ref{me1} it follows
that
\begin{equation}\label{pushforwardofEL}
(\psi_* \pi_{EL})((k,P), (k, Q)) = \langle \Omega_{(k_1,k_2)}(P^{k_1^{-1}},-P^{k_1^{-1}}), (Q^{k_1^{-1}},-Q^{k_1^{-1}})\rangle_{\u}
\end{equation}
where $\Omega_{(k_1,k_2)} (X,-X)=\{\Ad(k_1^{-1},k_2^{-1})\circ \H_{\k\times\k} \circ \Ad(k_1,k_2) (X,-X)\}_{i\p}$.
The result of applying the transformation $\Omega_{(k_1,k_2)}$ to
$(P^{k_1^{-1}},-P^{k_1^{-1}})$ is equivalent to the $i\p$ part of
\begin{equation}\label{onestepcloser}
\Ad(k_1^{-1},k_2^{-1})\circ \H_{\k\times\k} (P,-P^{k^{-1}}) =
((\H_\k (P))^{k_1^{-1}} , - (\H_\k (P^{k^{-1}}))^{k_2^{-1}}).
\end{equation}
Note that the left and right hand side of (\ref{onestepcloser}) are in $\u=\k\times\k$.
Temporarily denote the right hand side of (\ref{onestepcloser}) by $Z$.
Since $Z\in \u$, the projection to $i\p$ is given by $\frac{1}{2}(Z-Z^\theta)$ or
\begin{equation}\label{ippart}
 \frac{1}{2} \left( (\H_\k(P))^{k_1^{-1}}+(\H_\k (P^{k^-1}))^{k_2^{-1}} , - ((\H_\k(P))^{k_1^{-1}}+(\H_\k (P^{k^-1}))^{k_2^{-1}})\right).
\end{equation}
Substituting (\ref{ippart}) into the right hand side of (\ref{pushforwardofEL})
yields the expression
\begin{equation*}
\langle (\H_\k(P))^{k_1^{-1}}+(\H_\k (P^{k^-1}))^{k_2^{-1}} , Q^{k_1^{-1}} \rangle
\end{equation*}
from which it follows that
\begin{equation*}
(\psi_*\pi_{EL})((k,P),(k,Q)) = \langle (\H_\k + \Ad(k)\circ \H_\k \circ \Ad(k^{-1}))(P), Q  \rangle .
\end{equation*}
completing the proof of the theorem.
\end{proof}

For a compact group $K$ there are essentially two
Poisson structures intimately
related to the Lie theory of $K$.  On the one hand there is the Lu-Weinstein
Poisson Lie group structure whose symplectic foliation respects the
Bruhat decomposition of $K$.  This is given by the difference of the left
and right invariant bivector fields generated by $\H_\R$. On the other hand there is the Evens-Lu homogeneous
Poisson structure on $K$ whose symplectic foliation respects the
Birkhoff decomposition of $K$.  This is given by the sum of the
left and right invariant bivector fields generated by $\H_\R$.

Returning to the general case where $X$ is not necessarily a group, note that
all stability subgroups of $U$ corresponding to points in $X$ are
conjugate in $U$. Thus, the class of $\theta$ in the outer
automorphism group of $U$ is an invariant of $X$.  When this class is
trivial, the Evens-Lu Poisson structure is non-degenerate on an open
dense subset of points of $X$ as will be shown below.  This stems from
the fact that when $\theta$
is an inner automorphism, each Cartan subalgebra of $\k$ is, in fact,
a Cartan subalgebra of $\u$.  Thus $\h\cap \u=\t$ and $\h\cap \k=\t_0$
are equal.  Many statements simplify dramatically in the inner case.  For
example, the torus $T_w=\exp(\ker\{\Ad(w)\theta|_t-1\})$ which acts on the
layer of the Birkhoff decomposition of $X$ corresponding to $w$ admits
a much simpler description in the inner case as $T_0\cap (T_0)^w$.

\begin{thm}\label{innercase_symplectic_leaves}
In symmetric spaces for which $\theta$ is an inner automorphism, each
connected component of the layer of the Birkhoff decomposition of $X$
corresponding to the trivial element of the Weyl group is an open symplectic
leaf.  The components are indexed by the elements of order two in $T_0$.
\end{thm}
\begin{proof}
From corollary \ref{corollary_on_symplectic_leaves}:
the symplectic leaves of maximal dimension foliate the connected components
of the layer corresponding to $w=T_0\in W=N_U(T_0)/T_0$.  For brevity, this
layer will be referred to as the \emph{top layer}.  The leaves are
indexed by the elements $\w\in \phi(U/K)\cap T_0$.  Such an element
satisfies the equation $\w^{-1}=\w^\theta$, but $T_0$ is fixed by
$\theta$, so $\w$ must be an element of order 2.  Theorem 3 in \cite{pickrell}
shows that for each such $\w$ there exists an element $\w_1\in N_U(T_0)$
such that $\phi(\w_1K)=\w$.

In theorem \ref{me1}
it was shown that with the presentation of $T^*(U/K)\simeq U\times_K i\p$
the Evens-Lu Poisson bivector can be expressed relative to the Killing form
as
\begin{equation*} \label{bivector2}
\pi([u,X],[u,Y])=\langle \Omega_u X, Y\rangle
\end{equation*}
where $X,Y\in i\p$ and $\Omega_u\colon i\p \to i\p$ is given by
\begin{equation*}
\Omega_u(X)=\{\Ad(u^{-1})\circ \H\circ \Ad(u) (X)\}_{i\p}.
\end{equation*}
The kernel of $\Omega_{\w_1}(X)$ is equal to $\t_0^{\w_1^{-1}}\cap i\p
=\t_0 \cap i\p=0$.  Thus, the leaves in $X$ whose $\phi$-images lie in the
top layer of the Birkhoff decomposition are open.
\end{proof}

The elements of order two in $T_0$ are precisely the preferred basepoints
in $X$ whose existence was established in lemma \ref{basepoint_choice}.

From the classification of symmetric spaces in \cite{Helgason}, the list
of irreducible compact symmetric spaces for which the involution is an inner
automorphism includes, but is not limited to, the compact
Hermitian symmetric spaces.  Complex variables will be used in the
following section to exhibit locally
the Evens-Lu Poisson structure on some spaces of this type.

\section{Examples\label{examples}}

In this section, local expressions are recorded for the
Evens-Lu Poisson structure in a number of explicit examples.  The first
example is the complex Grassmannian, i.e. the space
$m$-planes in in $\C^{m+n}$.  As a symmetric space this may be presented as the quotient of the compact group $U=\SU(m+n)$
by the closed subgroup $K=\mathrm{S}(\mathrm{U}(m)\times \mathrm{U}(n))$.  This presentation arises from the natural
action of $\SU(m+n)$ on
$\C^{m+n}$ by linear isometries of the standard Hermitian inner product.  This descends to a
transitive action on the set of
complex $m$-planes through the origin.

The complexification of $\u=\su(m+n)$ is $\g=\sl(m+n,\C)$.  Let $\h$
be the diagonal matrices in $\sl(m+n,\C)$
and the triangular decomposition be the usual one where $\n_+$ consists
of the strictly upper triangular matrices
and $\n_-$ the strictly lower triangular matrices.  This $\u$-compatible
triangular decomposition generates
a standard Poisson Lie group structure on $U$.  The corresponding
Birkhoff decomposition of $G=\SL(m+n,\C)$
corresponds to the factorization produced in linear algebra through
Gaussian elimination.

The point in the Grassmannian corresponding to the plane spanned by the
first $m$ standard basis vectors
is an example of a preferred basepoint.  Denote this plane by $\C^m$.
The stability subgroup of this basepoint consists of the special unitary
transformations preserving this plane (and by necessity its Hermitian orthogonal
complement), i.e. $\mathrm{S}(\mathrm{U}(m)\times \mathrm{U}(n))$.
From this presentation one can readily see that
$\theta$ is the automorphism which negates the off-diagonal blocks.
This involution is an inner automorphism, given by conjugation by
a scalar multiple of the block diagonal matrix with an
$m\times m$ identity matrix in the upper diagonal block and an
$n\times n$ diagonal matrix with negative ones on the diagonal in the lower block.
When $n$ is even, the scalar multiple is one.
When $n$ is odd, the multiple is a primitive $2(m+n)$-th root of unity.

The complex Grassmannian is, additionally, a Hermitian symmetric space.
It is diffeomorphic to the quotient of $\SL(m+n,\C)$ by the parabolic
subgroup of the upper block triangular matrices of the form
\begin{equation*}
\left(\begin{array}{cc}
 A & B \\
 0 & D \\
\end{array}\right)
\end{equation*}
where $A$ is $m\times m$ and $D$ is $n\times n$.  It is through an
identification such as this that the Grassmannian inherits a complex
structure. Holomorphic coordinates can thus be used to present local
formulas.

The graph of a $\C$-linear transformation $Z\in \L(\C^m,\C^n)$ in
$\C^{m+n}$ written
\begin{equation*}
\{(X,ZX)\colon X\in \C^m\}
\end{equation*}
uniquely determines a point in the Grassmannian.  In fact, every complex $m$-dimensional subspace of $\C^{m+n}$ which is
transverse to $(\C^m)^\perp=\C^n$ can be realized in this way.
In this fashion $\L(\C^m,\C^n)$ provides an affine coordinate chart for
the Grassmannian with each point in a Zariski open subset of the
Grassmannian described by an $n\times m$ matrix $Z$ of complex numbers.

   Each coset of $\mathrm{S}(\U(n)\times \U(m))$ corresponding the graph
of a linear transformation contains a unique element $u$ which
has positive definite diagonal blocks.  This can be seen as follows.  Apply polar decomposition
to the diagonal blocks $A$ and $D$ of a special unitary matrix $u$, writing
$A=|A|P_A$ and $D=|D|P_D$ and then factor $u$ to find
\begin{equation}
u=\left(\begin{array}{cc} A & B \\ C & D \end{array}\right) =  \left(\begin{array}{cc} |A| & BP_D^* \\ CP_A^* & |D|\end{array}\right)
\left(\begin{array}{cc}P_A & 0 \\ 0 & P_D \end{array}\right) \label{factorization}
\end{equation}
The diagonal blocks of a special unitary matrix corresponding to the graph of a linear transformation are invertible,
so their polar decomposition is unique.   Furthermore, the diagonal
blocks of a special unitary matrix have conjugate determinants
(see proposition 1.3 of \cite{Forman}),
thus $P_A$ and $P_D$ have conjugate determinants.  Therefore, the
factorization in (\ref{factorization}) produces a unique
representative for $u$ modulo $\mathrm{S}(\mathrm{U}(m)\times \mathrm{U}(n))$ with the desired properties.

  This preferred coset representative can then be determined as a function of $Z$.  It is given by
\begin{equation*}
u(Z)=\left(\begin{array}{cc}(1+Z^*Z)^{-1/2} & -(1+Z^*Z)^{-1/2}Z^*  \\
Z(1+Z^*Z)^{-1/2} & (1+ZZ^*)^{-1/2} \end{array}\right).
\end{equation*}
This formula for $u(Z)$ is will be referred to as the
\emph{canonical representative} for the coset corresponding
to $uK$ depending on $Z$.

Cotangent vectors at $Z$ will be represented by $m\times n$ complex
matrices using the identification of the real cotangent space to the
Grassmannian with the holomorphic cotangent space which, at a point, is
further identified with $\L(\C^n,\C^m)$.

\subsection{The Grassmannian}
With respect to these coordinates the action of the bivector $\pi$ on
two cotangent vectors represented by
complex $m\times n$ matrices $V$ and $W$ may be computed by
\begin{equation*}
\pi (V,W)=i[\tr(L_ZV)^*W)-\tr((L_ZV)W^*)] \label{grassmann_formula}
\end{equation*}
where $L_Z$ is the $\R$-linear transformation $\L(\C^m,\C^n)\to \L(\C^m,\C^n)$ given by
\begin{eqnarray}
L_ZV & = & V-Z^*ZVZZ^* \nonumber \\
& & + Z^*((ZV-V^*Z^*)_+ + c.t.) \nonumber \\
& & - ((Z^*V^*-VZ)_+ + c.t ). \label{third_term}
\end{eqnarray}
In the above expression, $(\cdot)_+$ denotes the upper triangular part as before, $V^*$ denotes the conjugate
transpose of the matrix $V$, and $c.t.$ denotes the conjugate transpose of the preceding term.
This local formula is obtained by direct calculation from the equivariant
formula for $\pi$ in theorem \ref{me1}.

\subsection{Complex Projective Space}
   Complex projective space of dimension $n$, denoted $\CP^n$, is the space of complex lines through the origin in $\C^{n+1}$ and is
the Grassmannian with $m=1$.  In this case, the coordinate $Z$ is a column
vector and the matrices representing cotangent vectors are
row vectors.  The quantity $Z^*Z$ is a scalar which, for brevity, we
write as
\begin{equation*}
Z^*Z=\|Z\|^2=|z_1|^2+\dots +|z_n|^2.
\end{equation*}
Furthermore, the third term of $L_Z V$ in (\ref{third_term}) vanishes
as the matrix $Z^*V^*-VZ$ is one by one and thus has no
upper triangular part.  In this case the bivector can be
written more explicitly as
\begin{eqnarray*}
\pi & = & -i\left\{\sum_{j=1}^n S_j \delz{j}\wedge \delzbar{j} \right. \\
&  & \hspace{-3mm}+\left.\left(\sum_{j<k}z_jz_k\delz{j}\wedge\delz{k} - \sum_{j<k} z_j\zbar_k\|Z\|^2 \delz{j}\wedge \delzbar{k}\right) - c.c.\right\}
\end{eqnarray*}
where
\begin{equation*}
S_j = 1 +\sum_{k=1}^{j-1}|z_k|^2-|z_j|^2\|Z\|^2-\sum_{k=j+1}^n |z_k|^2.
\end{equation*}

  It is interesting to consider $\CP^2$ and $\CP^1$ in
further detail.  It appears that
the coefficients $S_j$ are reducible polynomials in the variables $|z_j|^2$ only
in these cases.  It is not
clear what the significance of this is (if any).
\begin{figure}[t]
\[\begin{xy}
\xyimport(144,144){\includegraphics{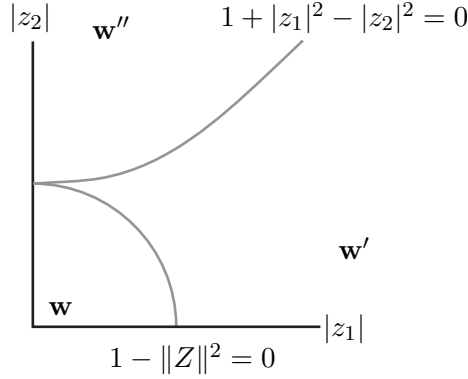}}
,(80,-15)*{1-\|Z\|^2=0}
,(155,155)*{1+|z_1|^2-|z_2|^2=0}
,(15,10)*{\w}
,(160,40)*{\w'}
,(40,150)*{\w''}
,(155,0)*{|z_1|}
,(0,155)*{|z_2|}
\end{xy}\]
\caption{The degeneracy locus for $\CP^2$.  The three elements $\w=\mathrm{diag}(1,1,1)$, $\w'=\mathrm{diag}(-1,-1,1)$,
$\w''=\mathrm{diag}(-1,1,-1)$ in $T_0^{(2)}\cap \phi(U/K)$ are depicted with
each component of the top stratum they determine.}
\end{figure}
For $\CP^2$ the bivector can be expressed locally as
\begin{eqnarray*}
\pi & = & -i\left\{S_1\delz{1}\wedge\delzbar{1}+ S_2\delz{2}\wedge\delzbar{2}\right. \\
&  & + \left.\left( z_1z_2\delz{1}\wedge\delz{2} - z_1\zbar_2\|Z\|^2 \delz{1}\wedge \delzbar{2}\right)-c.c.\right\}
\end{eqnarray*}
where $S_1=(1+|z_1|^2)(1-\|Z\|^2)$ and $S_2=(1-|z_2|^2)(1+\|Z\|^2)$.
  Since the symplectic leaves are open, one can write down the induced
symplectic structure in these coordinates by inverting the Poisson tensor.
Explicitly,
\begin{eqnarray*}
\omega & = & \frac{i}{p(Z,Z^*)}\left\{-S_2 dz_1\wedge d\zbar_1 -S_1 dz_2\wedge d\zbar_2 \right. \\
&  & + \left.\left(z_1z_2d\zbar_1\wedge d\zbar_2+z_1\zbar_2\|Z\|^2d\zbar_1\wedge dz_2\right)-c.c\right\}
\end{eqnarray*}
where $p(Z,Z^*)=(1+|z_1|^2-|z_2|^2)(1-\|Z\|^2)(1+\|Z\|^2)$.  The degeneracy locus of $\pi$ is given by the
variety $p(Z,Z^*)=0$.

To witness the connection with triangular factorization the
canonical representatives can be used to compute the complement of the
top layer of the Birkhoff decomposition
in these coordinates.  By introducing the real analytic function
$\varphi = (\sqrt{1+\|Z\|^2}-1)/\|Z\|^2$, one can compute the matrix
$(1+ZZ^*)^{-1/2}$.  Then
\begin{equation*}
u = \frac{1}{\sqrt{1+\|Z\|^2}} \left(
\begin{array}{ccc}1 & -\zbar_1 & -\zbar_2 \\
                  z_1 & 1+|z_2|^2\varphi & -z_1\zbar_2\varphi \\
                  z_2 & -\zbar_1 z_2\varphi & 1+|z_1|^2\varphi
\end{array}\right)
\end{equation*}
and
\begin{equation}\label{image_of_CP2}
uu^{-\theta} = \frac{1}{1+\|Z\|^2}\left(
\begin{array}{ccc}
1-\|Z\|^2 & -2\zbar_1 & -2\zbar_2 \\
2z_1 & 1-|z_1|^2+|z_2|^2 & -2z_1\zbar_2 \\
2z_2 & -2\zbar_1 z_2 & 1+|z_1|^2-|z_2|^2
\end{array}\right).
\end{equation}
The matrix $uu^{-\theta}$ in (\ref{image_of_CP2}) is in
$\phi(U/K)\cap \Sigma_1^G$ provided the principal minors are
non-vanishing.  Thus, the complement of the top layer is given in
coordinates by the vanishing locus of the product of the principal minors.
One can check that, in this case, this product is given by the smooth rational
function $p(Z,Z^*)/(1+\|Z\|^2)^3$.

In general, the degeneracy locus for $\pi$ on the complex Grassmannian
is given by the vanishing locus of a reducible
polynomial whose factors are given by the explicit formulas of Habermas (\cite{Habermas}).

\subsection{$\SU(2)$ and $\CP^1$}

For $X=\CP^1$ the group $U$ is $\SU(2)$ and $\g=\mathrm{sl}(2,\C)$.
As is typical with this example, everything can be computed.  The Lu-Weinstein
Poisson Lie group structure, the Evens-Lu Poisson structure that $\SU(2)$ inherits
as a symmetric space, and the Evens-Lu Poisson structure on $\CP^1$ will
all be displayed.  The end of this subsection returns to a topic discussed in the
introduction, relating the structure produced by the Evens-Lu structure
to the one produced by the Foth-Lu construction (\cite{FL}).  To conclude,
the Evens-Lu Poisson structure is exhibited as an element of the one
parameter family from \cite{KKR}.

Consider the standard triangular
decomposition of $\g$ as for the Grassmannian: $\h$ is
the set of traceless diagonal matrices and $\n_\pm$ are
the strictly upper triangular (resp. lower triangular) matrices.
Temporarily, set
\begin{equation*}
E_-=\left(\begin{array}{cc}0 & 0 \\ 1 & 0 \end{array}\right),\,\,\,\,H=\left(\begin{array}{cc}i & 0 \\ 0 & -i \end{array}\right), \text{ and }E_+=\left(\begin{array}{cc}0 & 1 \\ 0 & 0 \end{array}\right).
\end{equation*}
Then $\h=\span_\C\{H\}$, $\n_\pm=\span_\C\{E_\pm\}$, and $\u=\span_\R\{H,X,Y\}$ where $X=E_+-E_-$,
and $Y=i(E_++E_-)$. The triangular decomposition is stable with respect to the involution
selecting the stability subgroup $\mathrm{S}(\mathrm{U}(1)\times \mathrm{U}(1))$.

Using right translation, the tangent bundle to $U$ can be identified
with $U\times \u$.  Thus a bivector field
on $U$ can be identified with a smooth map
$U\to \bigwedge^2\u$.  In this presentation the value of the Lu-Weinstein
Poisson Lie group structure $\pi_U$ at an element
\begin{equation*}
\left(
\begin{array}{cc}
a & b \\
-\bar{b} & \bar{a}
\end{array}
\right)\in \SU(2)
\end{equation*}
is given by
\begin{equation}
\pi_U = (1-|a|^4+|b|^4)X\wedge Y + 2\Im(\bar{a}b)Y\wedge H - 2\Re(a\bar{b})H\wedge X  \label{PLstructure}
\end{equation}
whereas the Evens-Lu Poisson structure that $\SU(2)$ inherits as a symmetric
space is given by
\begin{equation*}
\pi^{\text{EL}} = (1+|a|^4-|b|^4)X\wedge Y + 2\Im (ab)Y\wedge H - 2\Re (ab) H\wedge X. \label{ELstructureSU(2)}
\end{equation*}
One can see that the Evens-Lu Poisson structure $\pi^{EL}$ vanishes precisely
when the principal minor, i.e., $a$, is zero.  This shows that the symplectic
foliation of $\pi^{EL}$ respects the Birkhoff decomposition of $\SU(2)\subset \SL(2,\C)$.

With the presentation of $\CP^1$ as $\SU(2)/\mathrm{S}(\U(1)\times \U(1))$ the
Evens-Lu construction gives the Poisson structure
\begin{equation}
\pi = -i(1-|z|^4)\delz{}\wedge\delzbar{}.\label{ELforCP1}
\end{equation}
The top layer of the Birkhoff decomposition has two connected components,
the upper and lower hemispheres.  The degeneracy locus of $\pi$ is the
equator.

\begin{figure}[h]\label{strataimage}
\[\begin{xy}
\xyimport(216,144){\includegraphics{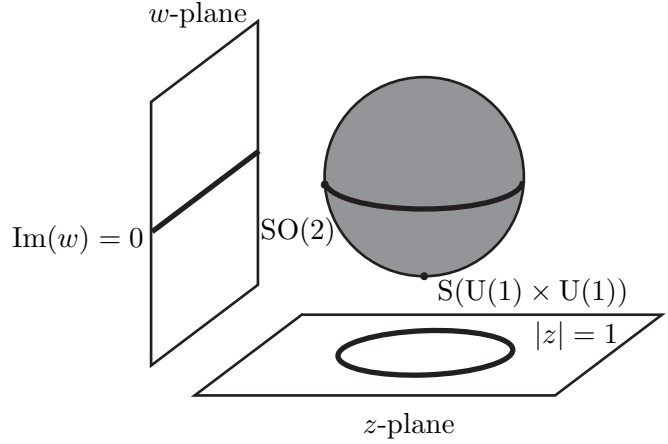}}
,(20,145)*{w\text{-plane}}
,(108,-10)*{z\text{-plane}}
,(62,63)*{\SO(2)}
,(160,40)*{\mathrm{S}(\mathrm{U}(1)\times \mathrm{U}(1))}
,(-30,60)*{\Im (w)=0}
,(178,24)*{|z|=1}
\end{xy}\]
\caption{The symplectic foliation for the $(U,\pi_U)$-homogeneous Poisson structure on the $\CP^1$.
The degeneracy locus is pictured in the $z$-coordinate when the stability subgroup of the basepoint is
$\mathrm{S}(\mathrm{U}(1)\times \mathrm{U}(1))$ and the $w$-coordinate when the stability subgroup is $\SO(2)$.
}
\end{figure}

Alternatively, Foth and Lu choose the basepoint of this symmetric space so that the stability subgroup is
$\SO(2)\subset \SU(2)$.  In this case, the corresponding non-compact real
form $\g_0$ is $\sl(2,\R)$ and
the Borel subalgebra $\mathfrak{b}=\h+\n_+$ is Iwasawa relative to $\g_0$.
With this choice of basepoint the quotient map $U\to U/K$ projects $\pi_U$
to a $(U,\pi_U)$-homogeneous Poisson structure on $U/K$.  The map
\begin{equation*}
w\mapsto \frac{1}{\sqrt{1+|w|^2}}\left(
\begin{array}{cc}
1+i\Re(w) & i\Im(w) \\
i\Im(w) & 1-i\Re(w)
\end{array}
\right)
\end{equation*}
gives a local cross section of the projection $\SU(2)\to \SU(2)/\SO(2)$.
In this coordinate, the projection of $\pi_U$ (\ref{PLstructure}) is given by
\begin{equation}
\pi = -2i\Im(w)(1+|w|^2)\delw \wedge \delwbar.
\end{equation}
These two points of view are illustrated in figure 3.  There, one can see
that the Foth-Lu construction and the Evens-Lu
construction produce the same $(U,\pi_U)$-homogeneous Poisson structure
on $\CP^1$.

In \cite{KKR} the existence of a parabolic subgroup of $G$ such that
$U\cap P$ is a Poisson Lie subgroup was established under the assumption
that $X$ was a Hermitian symmetric space.  The projection of $\pi_U$
under the quotient map
$U\mapsto U/(U\cap P)$ then defines a $(U,\pi_U)$-homogeneous Poisson
structure on $X$ which is compatible with the invariant Poisson
structure $\pi_{KKS}$.

In this example, the parabolic subgroup $P$ is actually the Borel subgroup
$HN^+\subset \SL(2,\C)$ and $\SU(2)\cap P$ is the diagonal torus
$\mathrm{S}(\U(1)\times \U(1))$.
Pushing forward $\pi_U$ from (\ref{PLstructure}) under the quotient map
gives the Poisson structure
\begin{equation}
\pi_{PL} = 2i|z|^2(1+|z|^2)\delz{} \wedge \delzbar{}
\end{equation}
which is degenerate only at the basepoint.
Furthermore, the invariant Poisson structure of Kostant-Kirillov-Souriau
is given in the $z$-coordinate by
\begin{equation}
\pi_{KKS} = i(1+|z|^2)^2 \delz{}\wedge \delzbar{}.
\end{equation}
The reader is invited to check that the Evens-Lu Poisson structure $\pi$
in (\ref{ELforCP1}) is equal to $\pi_{PL}+\lambda \pi_{KKS}$ with $\lambda=-1$.

\end{document}